\documentclass[a4paper,12pt]{amsart}
\usepackage{fancyheadings,amsmath,amscd,amsthm,amsfonts,latexsym,mathrsfs}
\DeclareMathOperator{\RicM}{^{\it M}\!Ricci}

\DeclareMathOperator{\dif}{d}

\renewcommand{\H}{\mathcal{H}}
\newcommand{\V}{\mathcal{V}}

\def \G{\Gamma}
\def \l{\lambda}

\def \phi{\varphi}
\def \S{\Sigma}

\def \Re{\mathbb{R}}

\def \Co{\mathbb{C}\,}

\begin{document}

\title{Topological restrictions for circle actions\\ 
	and harmonic morphisms}
\author{Radu Pantilie and John C. Wood}
\thanks{The first author gratefully acknowledges the support of the O.R.S. 
Scheme Awards, the School of Mathematics of the University of Leeds, 
the Tetley and Lupton Scholarships and  the Edward Boyle Bursary.}
\address{University of Leeds\\Dept. of Pure Mathematics\\
Leeds LS2 9JT\\UK.}
\email{pmtrp@amsta.leeds.ac.uk}
\address{University of Leeds\\Dept. of Pure Mathematics\\
Leeds LS2 9JT\\UK.}
\email{j.c.wood@leeds.ac.uk}

\newtheorem{thm}{Theorem}[section]
\newtheorem{lem}[thm]{Lemma}
\newtheorem{cor}[thm]{Corollary}
\newtheorem{prop}[thm]{Proposition}

\theoremstyle{definition}
\newtheorem{defn}[thm]{Definition}
\newtheorem{rem}[thm]{Remark}
\newtheorem{exm}[thm]{Example}

\numberwithin{equation}{section}

\maketitle
\thispagestyle{empty} 

\section*{Abstract}
\begin{quote}
{\footnotesize
Let $M^m$ be a compact oriented smooth manifold which admits a smooth circle action 
with isolated fixed points which are isolated as singularities as well. Then all the 
Pontryagin numbers of $M^m$ are zero and its 
Euler number is nonnegative and even. In particular, $M^m$ has signature zero.\\ 
Since a non-constant harmonic morphism with one-dimensional fibres gives rise to a circle 
action we have the following applications:\\
(i)  many compact manifolds, for example $\Co\!P^{n}$, $K3$ surfaces, $S^{2n}\times P_g$ ($n\geq2$) 
where $P_g$ is the closed surface of genus $g\geq2$ can never be the total space of a non-constant 
harmonic morphism with one-dimensional fibres whatever metrics we put on them;\\
(ii) let $(M^4,g)$ be a compact orientable four-manifold and 
$\phi:(M^4,g)\to(N^3,h)$ a non-constant harmonic morphism. Suppose that one of the following 
assertions holds:\\
$\bullet$ $(M^4,g)$ is half-conformally flat and its scalar curvature is zero,\\ 
$\bullet$ $(M^4,g)$ is Einstein and half-conformally flat,\\ 
$\bullet$ $(M^4,g,J)$ is Hermitian-Einstein.\\ 
Then, up to homotheties and Riemannian 
coverings, $\phi$ is the canonical projection $T^4\to T^3$ between flat tori.}
\end{quote}

\section*{Introduction} 

\indent
It is well known that if a compact oriented smooth manifold $M$ admits a smooth \emph{free} circle 
action then its Euler number and all its Pontryagin numbers are zero. This follows from the 
fact that the tangent bundle of $M$ is the Whitney sum of a trivial real line bundle and the 
pull back of the tangent bundle of the orbit space.\\ 
\indent 
In this paper we generalize this by proving that 
if $M$ is a compact oriented smooth manifold which admits a smooth circle action with 
isolated fixed points which are isolated as singularities as well then (i) all the Pontryagin numbers 
of $M$ are zero (in particular, the signature of $M$ is zero), (ii) the Euler number of $M$ is even 
and is equal to the number of fixed points (Theorem \ref{thm:circle}).  We obtain this by using 
a well known formula of R.~Bott \cite{Bott} (see also \cite{Kob}). Also, we apply an idea 
of J.D.S.~Jones to prove (Theorem \ref{thm:top4d}) that the signature of a compact 
oriented $4$-manifold endowed with a non-trivial circle action for which each fixed point 
has equal exponents is equal to the Euler number of the normal bundle of the components of 
dimension $2$ of the fixed point set.

\indent 
\emph{Harmonic morphisms} are 
smooth maps $\phi:(M,g)\to(N,h)$ between Riemannian manifolds which preserve Laplace's equation. 
They are characterised 
as harmonic maps which are \emph{horizontally weakly conformal} \cite{Fug}, 
\cite{Ish}, 
i.e.\ for each point $x\in M$ either $\dif\!\phi_x=0$ or $\dif\!\phi_x$ is 
surjective and  
maps $\H_x=({\rm ker}\dif\!\phi_x)^{\perp}$ conformally onto $T_{\phi(x)}N$. 
Classification 
results for harmonic morphisms with one-dimensional fibres appear in \cite{BaiWoo-Bernstein}, 
\cite{BaiWoo-spfo}, \cite{BaiWoo1}, \cite{Bry}, \cite{Pan}, \cite{Pan-4to3}. 
In \cite{Bai} it is proved that any non-constant harmonic morphism with one-dimensional fibres 
defined on a Riemannian manifold of dimension at least five is submersive whilst, for domains 
of dimension four, only isolated critical points can occur. Moreover, in \cite{Bai} it is proved 
that any non-constant harmonic morphism $(M^4,g)\to(N^3,h)$ induces a locally smooth circle 
action on $M^4$. We show in fact that, at least outside the critical set, 
the action is smooth and free (consequence of Proposition \ref{prop:fact}). 
It follows that the 
result of Theorem \ref{thm:circle} can be applied to obtain topological restrictions for 
the total space of a harmonic morphism with one-dimensional fibres. These are obtained in 
Theorem \ref{thm:restr} from which it immediately follows that  $\Co\!P^n$, $K3$ surfaces, 
$S^{2n}\times P_g$ 
(where $P_g$ is the closed (real) surface of genus $g$ and $n\geq2$) can 
never be the total space of a non-constant harmonic morphism with 
one-dimensional fibres whatever metrics we put on them. The result regarding $\Co\!P^2$ 
answers a question formulated by P.~Baird in a conversation.\\ 
\indent
Further applications are obtained in Theorem \ref{thm:EinsteinASDscalar-flat} where we show 
that, up to homotheties and Riemannian coverings, the canonical projection $T^4\to T^3$ between flat tori is 
the only harmonic morphism with one-dimensional fibres which is defined on a compact 
half-conformally flat $4$-manifold which is either Einstein or scalar-flat.  We then 
apply this result to obtain a new proof of a result of \cite[Theorem 4.11]{Pan-4to3} 
in which the same conclusion is proved for a harmonic morphism $\phi:(M^4,g)\to(N^3,h)$ 
between compact Einstein manifolds. To obtain this new proof we first show 
(Proposition \ref{prop:ASD}) 
that if $\phi:(M^4,g)\to(N^3,h)$ is a non-constant harmonic morphism between orientable Einstein 
manifolds then $(M^4,g)$ is half-conformally flat.\\ 
\indent
In Theorem \ref{thm:Hermitian-Einstein}, we prove that the same conclusion as that in 
Theorem \ref{thm:EinsteinASDscalar-flat} holds for a harmonic morphism 
$\phi:(M^4,g,J)\to(N^3,h)$ defined on a compact Hermitian-Einstein four-manifold. 

\indent
We are grateful to A.L.~Edmonds, I.~Hambleton, J.D.S.~Jones and M.~McCooey for useful information on circle actions. We are also grateful to P.~Baird for 
pointing out a mistake in Proposition \ref{prop:fact}.

\section{A topological restriction for circle actions}

\indent
By a \emph{singularity} of a group action we mean a point at which the isotropy group is 
non-trivial. A \emph{fixed point} is a singularity where the isotropy group is the entire 
group.\\ 
\indent
Let $M^m$ be a compact oriented smooth manifold of dimension $m\geq1$ endowed with a smooth 
circle action. Let $F$ denote its fixed point set and $V$ its infinitesimal generator. 
Let $g$ be a Riemannian metric on $M$ with respect to which $V$ is a Killing vector field. 
Such a metric can be obtained by averaging an arbitrary Riemannian metric over the action.  Let $\nabla$ denote the Levi-Civita connection on $(M^m,g)$.\\ 
\indent 
Obviously $F$ is the zero set 
of $V$ and thus, its connected components are totally-geodesic submanifolds of $(M,g)$ 
of even codimension (see \cite{Kob}).\\  
\indent
Let $x\in F$ and suppose that the connected component $N$ of $x$ in $F$ has codimension 
$2r$. Because $(\nabla V)_x$ is a skew-symmetric endomorphism of 
$(T_xM,g_x)$, with respect to a suitably chosen orthonormal frame, $(\nabla V)_x$ 
is represented by the direct sum of the zero square matrix of dimension $m-2r$  
and $\oplus_{j=1}^r \bigl(\begin{smallmatrix}0&-m_j\\m_j&0\end{smallmatrix}\bigr)$ where 
$m_j>0$. In fact, from \cite[Chapter I, Proposition 1.9]{KoNo} it follows that  $m_j\in\mathbb{Z}$  
since $V$ integrates to give an $S^1$ action on $M^m$.  
Indeed, since $V$ is Killing its flow commutes with the exponential map; hence, via the exponential map 
at $x$, the linear flow induced by $(\nabla V)_x$ on $T_xM$ is locally equivalent to the flow of $V$. 
Following \cite{JonRaw} we shall call the (positive) integers $m_j$ the \emph{exponents} of the action 
at the fixed point  $x$. This is, of course, motivated 
by the fact that the exponential map of $(M^m,g)$ at $x$ induces a local equivalence 
between the given $S^1$ action and the following $S^1$ action on 
$\Re^m=\Re^{m-2r}\oplus{\Co^r}$: 
$$z\cdot(x_1,\ldots,x_{m-2r},z_1,\ldots,z_r)=(x_1,\ldots,x_{m-2r},z^{m_1}\,z_1,\ldots,z^{m_r}\,z_r)\;.$$ 
(In particular, this shows that the exponents $(m_1,\dots,m_r)$ do not depend on the metric $g$.)  
If $x\in F$ is an isolated fixed point (equivalently $m=2r$) then the orientation induced by the 
corresponding orthonormal frame is determined by the $r$-tuple $(m_1,\ldots,m_r)$. 
Let $\epsilon(x)$ be $+1$ or $-1$ according to whether or not 
this orientation agrees with the orientation of $T_xM$ (cf.\ \cite{JonRaw}).\\

\indent
The main result of this section is the the following. 

\begin{thm} \label{thm:circle} 
Let $M^m$ be a compact oriented smooth manifold which admits a smooth circle action whose 
fixed points are isolated singularities.\\ 
\indent
Then {\rm (i)} all the Pontryagin numbers of $M^m$ are zero, {\rm (ii)} the Euler number of $M^m$ is even 
and is equal to the number of fixed points. In particular, the signature of $M^m$ is zero. 
\end{thm} 
\begin{proof} 
\indent
Let $x\in F$. Because $x$ is an isolated singularity the exponents at $x$ are all equal to $1$. 
Equivalently, there exists an orthonormal basis of $(T_xM,g_x)$ with 
respect to which the matrix of $(\nabla V)_x$ is the direct sum of $n$ copies of 
$\bigl( \begin{smallmatrix} 0&-1\\ 1&0 \end{smallmatrix} \bigr)$.\\  
\indent 
Thus, if $F\neq\emptyset$, then $\dim M=2n$ is even. Let $f$ be an 
${\rm Ad}({\rm SO}(2n))$-invariant symmetric polynomial of degree  $p\leq n$. 
Then, by a result of R.~Bott \cite[Theorem 2]{Bott} (see also \cite[Theorem II.6.1]{Kob})   
we have 
\begin{equation} \label{e:Bott} 
\sum_{x\in F}\frac{f((\nabla V)_x)}{\chi_n(\frac{1}{2\pi}(\nabla V)_x)}=\int_Mf(R) 
\end{equation} 
where $\chi_n$ is the Pfaffian (see \cite[p.\ 68]{Kob} or 
\cite[p.\ 309]{MilSta}) and $f(R)$ is the closed $2p$-form on $M$ which represents the 
cohomology class induced by the Chern-Weil morphism applied to $f$ via the Levi-Civita 
connection of $(M^{2n},g)$ (see \cite[Chapter XII]{KoNo} or \cite[Appendix C]{MilSta}). 
Note that the right hand side of \eqref{e:Bott} is zero if $p<n$.\\ 
\indent
It is easy to prove that 
\begin{equation} \label{e:easy} 
\chi_n(\tfrac{1}{2\pi}(\nabla V)_x)=\frac{(-1)^n\,\epsilon(x)}{(2\pi)^n}\;. 
\end{equation} 
\indent
By taking $f=1$, from \eqref{e:Bott} and \eqref{e:easy}, we obtain 
\indent
\begin{equation} \label{e:sumzero} 
\sum_{x\in F}\epsilon(x)=0\;. 
\end{equation}   
\indent
By taking $f=\frac{1}{(2\pi)^n}\,\chi_n$\,, from  \eqref{e:Bott} and 
the Gauss-Bonnet Theorem we obtain that the Euler number of $M^{2n}$ is equal to the cardinal 
of $F$ (this also follows from the Poincar\'e-Hopf theorem or from  
\cite[Theorem II.5.5]{Kob}). But, by \eqref{e:sumzero}, 
the cardinal of $F$ must be even and hence the Euler number of $M^{2n}$ is even.\\ 
\indent
By definition, if $\dim M$ is not divisible by four then all the Pontryagin numbers 
of $M$ are zero.\\ 
\indent
Suppose that $n=2p$ and let $i_1,\ldots,i_r$ be a partition of $p$. Denote by 
$p_{i_k}$ the ${\rm Ad}({\rm SO}(2m))$-invariant symmetric polynomial of degree $2i_k$ such that 
$p_{i_k}(\tfrac{1}{2\pi}\,R)$ represents the $i_k$'th Pontryagin class of $M$.\\ 
\indent
Let $x\in F$ and recall that $(\nabla V)_x$ is the direct sum of $n$ copies of 
$\bigl( \begin{smallmatrix} 0&-1\\ 1&0 \end{smallmatrix} \bigr)$. Then it is obvious 
that for any ${\rm Ad}({\rm SO}(2n))$-invariant 
symmetric polynomial $f$ we have 
\begin{equation} \label{e:Pontryagin} 
f((\nabla V)_x)=c(f,n) 
\end{equation} 
where $c(f,n)$ is a constant which depends just on $f$ and $n$ but not on $x\in F$.\\
\indent
By taking $f=p_{i_1}\!\cdots p_{i_r}$ in \eqref{e:Bott} it follows  from \eqref{e:easy}, 
\eqref{e:sumzero}, \eqref{e:Pontryagin} that all the Pontryagin numbers of $M^{2n}$ are zero. 
The fact that $M^{2n}$ has zero signature follows from Hirzebruch's Signature Theorem 
(see \cite{MilSta}).      
\end{proof}

\section{Blowing-up isolated fixed points of circle actions}

For $j=1,\ldots,r$ let $m_j$ be positive integers. We consider the following action 
on $\Co^r$: 
\begin{equation} \label{e:exmaction} 
z\cdot(z_1,\ldots,z_r)=(z^{m_1}\,z_1,z^{m_2}\,z_2,\ldots,z^{m_r}\,z_r)\;. 
\end{equation} 
Obviously, the fixed point set of this action is $\{0\}$.\\ 
\indent
In what follows we need the notion of equivariant connected sum of two manifolds 
endowed with circle actions. The idea of the following construction comes from  \cite{ChuLam}. 
\begin{defn} 
Let $M$ and $N$ be manifolds, $\dim M=\dim N=2r$, both endowed with non-trivial circle actions.\\ 
\indent
Let $x\in M$ and $y\in N$ be isolated fixed points of these actions having the same set of 
exponents $\{m_1,\ldots,m_r\}$. Further, assume that $\epsilon(x)=-\epsilon(y)$.\\ 
\indent
Then suitably chosen neighbourhoods $U_x$ and $U_y$ about $x$ and $y$ in $M$ and $N$, respectively, 
are equivariantly diffeomorphic to open balls of radius three about $0\in\epsilon(x)\Co^r$ and 
$0\in\epsilon(y)\Co^r$, respectively, where $\Co^r$ is endowed with the circle action of \eqref{e:exmaction}.\\ 
\indent 
The \emph{equivariant connected sum of $M$ and $N$ (about $x\in M$ and $y\in N$)} is the quotient 
 $(M\setminus V_x\sqcup N\setminus V_y)/_{\sim}$ where $V_x$ and $V_y$ 
correspond,  via the above equivariant diffeomorphisms, to the open balls $\epsilon(x)B(1)$ and 
$\epsilon(y)B(1)$, of radius $1$, whilst $\sim$ is induced by the identification 
$(2+t)u\sim (2-t)u$ with $u\in\partial B(1)$ and  $t\in[-1,1]$.\\ 
\indent
Thus the actions glue together to give a non-trivial circle action on the connected 
sum $M\# N$. 
\end{defn} 
\indent
\begin{exm} 
Let $\Co^r$ and $-\Co^r$ be endowed with the action given by \eqref{e:exmaction}. 
It is easy to see that the connected sum $\Co^r\#-\Co^r$, suitably constructed about $0\in\Co^r$ 
in each term, inherits in a canonical way a circle action (here $-\Co^r$ denotes $\Co^r$ considered 
with the orientation opposite to its usual one). Moreover this action is without fixed points. 
\end{exm}

\indent
The following definition is based on an idea of J.D.S.~Jones arisen from a private conversation. 
We formulate it only for isolated fixed points although it can be given for any connected component of the fixed point set. 

\begin{defn} 
Let $x\in F$ be an isolated fixed point with exponents $(m_1,\ldots,m_r)$. 
The \emph{blow-up} of $M$ (considered with the given action) at $x$ is the 
equivariant connected sum of $M$ and $-\epsilon(x)\Co\!P^r$, about $x\in M$ and $[1,0,\ldots,0]\in\Co\!P^r$, 
where $\Co\!P^r$ is considered with the $S^1$ action 
$$z\cdot[z_0,z_1,\ldots,z_r]=[z_0,z^{m_1}\,z_1,\ldots,z^{m_r}\,z_r]\;.$$ 
\end{defn} 

\indent
In what follows the following obvious lemma will play an important role. 

\begin{lem} \label{lem:codim2} 
Let $M$ be a manifold endowed with a circle action and let $F$ be its fixed point set.\\ 
\indent
Let $x$ be an isolated fixed point whose exponents are equal: $m_1=\ldots=m_r$.\\ 
\indent
Let $\widehat{F}$ be the fixed point set of the induced action on the blow-up of $M$ at $x$. 
Then $$\widehat{F}=(F\setminus\{x\})\cup\Co\!P^{r-1}$$ 
where $\Co\!P^{r-1}=\{\,[z_0,\ldots,z_r]\in\Co\!P^r\,|\,z_0=0\,\}$. 
\end{lem} 

\begin{rem} \label{rem:codim2} 
{}From the above lemma it follows that if besides isolated fixed points 
a circle action has only components of codimension two then after 
blowing up all the isolated fixed points we obtain a manifold endowed 
with a circle action whose fixed point set is of codimension two.\\ 
\indent
In particular, if the starting manifold is of dimension four then after 
blowing-up the isolated fixed points we obtain a manifold endowed with 
a circle action whose fixed point set is of dimension two. 
\end{rem} 

For the next lemma recall the 
${\rm Ad}({\rm SO}(4))$-invariant polynomial $p_1$ on $so(4)$ 
given by $p_1(A)=\sum_{i<j}(a^i_j)^2$ where $A=(a^i_j)_{i,j=1,\ldots,4}$. 
As is usual, we shall also denote by $p_1$ the corresponding 
${\rm Ad}({\rm SO}(4))$-invariant symmetric bilinear form on $so(4)$.\\   
\indent 
If $(M^4,g)$ is a Riemannian $4$-manifold then, by the Chern-Weil theorem, 
$p_1(\frac{1}{2\pi}R)$ represents the first Pontryagin class of $M^4$ 
where $R$ is the curvature form of the Levi-Civita connection of $(M^4,g)$ 
(see \cite{KoNo}, \cite{MilSta}).\\ 
\indent
Also we need the following definition (see \cite[p.\ 69]{Kob}). 

\begin{defn} \label{defn:residue} 
Let $(M^{2n},g)$ be a Riemannian manifold and let $V$ be a Killing vector field on it. 
Let $N$ be a component of codimension $2r$ of the zero set of $V$.\\ 
\indent
For an Ad(SO($2n$))-invariant polynomial $f$ of degree $n$ on $so(2n)$ the residue of $V$ 
over $N^{2n-2r}$ is given by 
\begin{equation} \label{e:residue} 
{\rm Res}_f(N)\,t^{n-r}=\int_N\frac{f(\frac{1}{2\pi}(tR+\nabla V))}
{\chi_r(\frac{1}{2\pi}(tR^{\perp}+(\nabla V)^{\perp}))}
\end{equation} 
where $\nabla$ is the Levi-Civita connection of $(M^{2n},g)$, $R$ its curvature form, 
$^{\perp}$ denotes the components in $End(TN^{\perp})$ and $\chi_r$ is the Pfaffian. 
(Here we expand the right hand side of \eqref{e:residue} as a power series in $t$.) 
\end{defn}

\begin{lem} \label{lem:restwo} 
Let $(M^4,g)$ be an oriented $4$-manifold and let $V$ be a Killing vector field 
on it. Suppose that the zero set of $V$ has a component $N^2$ of dimension 
two.\\ 
\indent
Then the residue ${\rm Res}_{p_1}(N)$ is given by 
$${\rm Res}_{p_1}(N)=\chi(TN^{\perp})[N]$$ 
where $\chi(TN^{\perp})[N]$ is the Euler number of the normal bundle $TN^{\perp}$ 
of $N$. 
\end{lem} 
\begin{proof} 
{}From Definition \ref{defn:residue} it follows that we can write 
\begin{equation} \label{e:resdef} 
{\rm Res}_{p_1}(N)\,t=\int_N\frac{p_1(\frac{1}{2\pi}(tR+\nabla V))}{\frac{1}{2\pi}(tR+\nabla V)^{\nu}}   
\end{equation} 
where 
$$\bigl(\begin{smallmatrix} x&y\\z&t \end{smallmatrix}\bigr)\oplus
\bigl(\begin{smallmatrix} 0&a\\-a&0 \end{smallmatrix}\bigr)^{\nu}=a\;,$$
and recall that $(\nabla V)|_N$ is a section of $End(TN)\oplus End(TN^{\perp})$ 
(see \cite[Chapter II, Theorem 5.3]{Kob}).  
In fact $(\nabla V)^{\nu}|_N:N\to\Re$ is the nowhere zero function on $N$ characterised by 
$$(\nabla V)|_N=\bigl(\begin{smallmatrix} 0&0\\0&0 \end{smallmatrix}\bigr)\oplus 
\bigl(\begin{smallmatrix}0&(\nabla V)^{\nu}\\-(\nabla V)^{\nu}&0\end{smallmatrix}\bigr)\;,$$  
then \eqref{e:resdef} becomes 
\begin{equation} \label{e:calc} 
\begin{split} 
{\rm Res}_{p_1}(N)\,t=&\frac{1}{2\pi}\int_N\bigl(t^2p_1(R,R)+2tp_1(R,\nabla V)+p_1(\nabla V,\nabla V)\bigr)\times\\ 
&\times\Bigl(\,\frac{1}{(\nabla V)^{\nu}}-t\,\frac{R^{\nu}}{((\nabla V)^{\nu})^2}+
t^2\frac{(R^{\nu})^2}{((\nabla V)^{\nu})^3}-\cdots\Bigr)\\ 
	=&\frac{t}{2\pi}\int_N\bigl(2p_1(R,U)-R^{\nu}\bigr) 
\end{split} 
\end{equation} 
where $U=\frac{1}{(\nabla V)^{\nu}}\nabla V|_N$ - note that with respect to a suitably chosen 
adapted orthonormal frame we have $U=\bigl(\begin{smallmatrix} 0&0\\0&0 \end{smallmatrix}\bigr)\oplus 
\bigl(\begin{smallmatrix}0&1\\-1&0\end{smallmatrix}\bigr)$. Also note that $J=U|_N$ is an 
almost complex structure on the normal bundle $TN^{\perp}$ (see \cite{Kob} for the general case).\\ 
\indent
Because $N$ is a totally-geodesic submanifold of $(M,g)$ we have that $p_1(R,U)=R^{\nu}$; hence,  
from \eqref{e:calc} it follows that 
\begin{equation} \label{e:restwo} 
{\rm Res}_{p_1}(N)=\frac{1}{2\pi}\int_NR^{\nu}\;.  
\end{equation} 
By the same reason we also have that 
$J\otimes R^{\nu}=\bigl(\begin{smallmatrix} 0&R^{\nu}\\-R^{\nu}&0 \end{smallmatrix} \bigr)$ 
is the curvature form of the connection induced by $\nabla$ on $TN^{\perp}$.\\ 
\indent
The proof now follows from the Chern-Weil theorem. 
\end{proof} 

\indent
We now state the main result of this section. 

\begin{thm} \label{thm:top4d} 
Let $M^4$ be a compact oriented $4$-manifold endowed with a non-trivial 
circle action.\\ 
\indent
Let $F=F_0\cup F_2$ be the fixed point set where $F_0$ is the set of 
isolated fixed points and $F_2$ is the union of the components of dimension $2$ 
of the fixed point set. Suppose that the exponents of each isolated fixed point $x\in F_0$ are equal.\\ 
\indent
Then the signature $\sigma[M]$ of $M$ is given by 
\begin{equation} \label{e:result} 
\sigma[M]=\sum_{x\in F_0}\epsilon(x)=\chi(TF_2^{\perp})[F_2]\;. 
\end{equation}  
\indent
In particular, the signature of $M$ is given by the Euler number of the normal bundle of the components of dimension 2 of the 
zero set of $V$. 
\end{thm} 
\begin{proof} 
Let $\widehat{M}$ be the manifold (endowed with a circle action)  obtained by blowing-up the isolated 
fixed points, i.e.\ the points of $F_0$.\\  
\indent
Then, since signatures add when taking connected sums, the signature $\sigma[\widehat{M}]$ of $\widehat{M}$ is given by 
\begin{equation} \label{e:signblow} 
\sigma[\widehat{M}]=\sigma[M]-\sum_{x\in F_0}\epsilon(x)\;. 
\end{equation} 
\indent
Because the exponents of each isolated fixed point are equal by hypothesis, the induced circle 
action 
on $\widehat{M}$ has \emph{no} isolated fixed points. This follows from Lemma \ref{lem:codim2} (see also 
Remark \ref{rem:codim2}).\\ 
\indent
Thus we can apply \cite[Theorem 4.2]{JonRaw} to obtain that $\widehat{M}$ has signature zero. Combining 
this with \eqref{e:signblow} gives 
\begin{equation} \label{e:firsteq} 
\sigma[M]=\sum_{x\in F_0}\epsilon(x)\;, 
\end{equation} 
i.e.\ the first equality of \eqref{e:result}.\\   
\indent
Now, take a metric on $M$ with respect to which $S^1$ acts by isometries. Then by applying the  Bott formula (see \cite[Chapter II, Theorem 6.1]{Kob}) to the infinitesimal generator of this action 
we obtain
\begin{equation} \label{e:p1blow} 
p_1[M]={\rm Res}_{p_1}(F_0)+{\rm Res}_{p_1}(F_2)\;. 
\end{equation} 
Since $p_1(A)=\pm2\chi_2(A)$ for $A\in so(4)_\pm$ (because $A=(a^i_j)\in so(4)_\pm$ 
if and only if $a^1_2=\pm a^3_4$, $a^1_3=\mp a^2_4$ and $a^1_4=\pm a^2_3$), 
we have ${\rm Res}_{p_1}(F_0)=2\sum_{x\in F_0}\epsilon(x)$.  
By using this fact and Lemma \ref{lem:restwo} the equation \eqref{e:p1blow} 
becomes 
\begin{equation} \label{e:inri} 
p_1[M]=2\sum_{x\in F_0}\epsilon(x)+\chi(TF_2^{\perp})[F_2]\;. 
\end{equation} 
\indent
But by Hirzebruch theorem $p_1[M]=3\sigma[M]$ which together with \eqref{e:firsteq} and \eqref{e:inri} 
gives 
\begin{equation} \label{e:q} 
3\sum_{x\in F_0}\epsilon(x)=2\sum_{x\in F_0}\epsilon(x)+\chi(TF_2^{\perp})[F_2] 
\end{equation} 
\indent
which immediately yields the second equality of \eqref{e:result}. 
\end{proof} 

\begin{rem} 
1) Obviously, Theorem \ref{thm:top4d} generalizes the result of \cite[Theorem 4.2]{JonRaw}. 
However, note that our proof of Theorem \ref{thm:top4d} uses that result.\\ 
\indent
Also, Theorem \ref{thm:top4d} generalizes the result of Theorem \ref{thm:circle} for 
dimension four. \\ 
\indent
2) See \cite{LaMi}, \cite{JonRaw} for some related results on circle actions. 
\end{rem}

\section{An application to harmonic morphisms} 

Let $\phi:(M^{n+1},g)\to(N^n,h)$\,, $n\geq1$\,, be a non-constant harmonic morphism between 
compact oriented Riemannian manifolds. Then, by a result of P.~Baird 
\cite{Bai} the set $\S$ of critical points of $\phi$ is empty if $\dim M\geq5$ and is 
discrete 
if $\dim M=4$\,. For $x\in M\setminus\S$ set $\V_x={\rm ker}\,\dif\!\phi_x$ and let 
$\H_x=\V_x^{\perp}$\,. 
The resulting distributions $\V$ and $\H$ on $M\setminus\S$ shall be called, as usual, the \emph{vertical} 
distribution and \emph{horizontal} distribution, respectively. 

\begin{prop} \label{prop:fact} 
Let $\phi:(M^{n+1},g)\to(N^n,h)\:\,(n\geq1)$ be a non-constant harmonic morphism between 
compact oriented Riemannian manifolds. Let $\S$ be the set of critical points of $\phi$\,.\\
\indent
Then $\phi|_{M\setminus\phi^{-1}(\phi(\S))}$ can be factorised as $\xi\circ\psi$ where $\psi:M\setminus\phi^{-1}(\phi(\S))\to P$ 
is the projection of an $S^1$ principal bundle (for which 
$\H|_{M\setminus\phi^{-1}(\phi(\S))}$ is a principal connection) 
and $\xi:P\to N\setminus\phi(\S)$ is a smooth covering 
projection. Moreover, for $n=3$, the smooth free $S^1$ action on 
$M\setminus\phi^{-1}(\phi(\S))$ extends to a continuous 
$S^1$ action on $M$ which is smooth over $M\setminus\S$ and whose fixed point set is $\S$\,.\\
\indent
Furthermore, if $k$ is the unique metric on $P$ with respect to which 
$\xi:(P,k)\to(N\setminus\phi(\S),h|_{N\setminus\phi(\S)})$ is a Riemannian 
covering then 
$\psi:(M\setminus\phi^{-1}(\phi(\S)),g|_{M\setminus\phi^{-1}(\phi(\S))})\to(P,k)$ is a submersive 
harmonic morphism with connected fibres. 
\end{prop}
\begin{proof} 
Because $M$ and $N$ are oriented, $\V$ is orientable. Thus we can choose $V\in\G(\V|_{M\setminus\S})$ such that $g(V,V)=\l^{2n-4}$\,. 
Clearly, $V$ is smooth on $M\setminus\S$. 
Furthermore, when $n=3$, since $|V|\to0$ as we approach a critical point, $V$ extends to a continuous vector field on $M$ whose zero set is $\S$ and  
the flow of $V$ extends to a continuous flow on $M$ whose fixed point set is $\S$\,.\\ 
\indent
For any $n\geq1$, it is easy to see that 
$\phi|_{M\setminus\phi^{-1}(\phi(\S))}$ is a proper submersion. Then by a  
well-known result of C.~Ehresmann \cite{Ehr}  
$\phi$ restricted to $M\setminus\phi^{-1}(\phi(\S))$ is the projection of a locally trivial fibre bundle. 
In particular, the orbit space $P$ of $V|_{M\setminus\phi^{-1}(\phi(\S))}$ 
is a smooth manifold. Thus, $\phi|_{M\setminus\phi^{-1}\phi(\S))}$ can be factorised as 
$\xi\circ\psi$ where  $\psi:M\setminus\phi^{-1}(\phi(\S))\to P$ has connected fibres and   
$\xi:P\to N\setminus\phi(\S)$ is a covering projection.\\ 
\indent
Let $k=\xi^*(h)$ be the unique metric on $P$ with respect to which 
$\xi:(P,k)\to(N\setminus\phi(\S),h|_{N\setminus\phi(\S)})$ becomes a Riemannian 
covering. It is obvious that 
$\psi:(M\setminus\phi^{-1}(\phi(\S)),g|_{M\setminus\phi^{-1}(\phi(\S))})\to(P,k)$ is a  submersive harmonic 
morphism with compact connected fibres. {}From \cite[Theorem 2.9]{Pan} it follows that $\psi$ 
is the projection of a circle bundle where the action on the total space $M\setminus\phi^{-1}(\phi(\S))$ is induced by the flow of $V|_{M\setminus\phi^{-1}(\phi(\S))}$.   
\end{proof} 

\begin{rem} 
Recall (\cite{Fug}, \cite{Ish}, \cite{BaiWoo2}) that, for $n=1$ a harmonic morphism from 
$(M^2,g)$ to $(N^1,h)$ is, essentially, a harmonic function.\\
\indent
When $n=2$ given a non-constant harmonic morphism $\phi:(M^3,g)\to(N^2,h)$ with $M^3$ compact 
the $S^1$ action extends \textit{smoothly} over the set of critical points and induces on $M^3$ a structure of a smooth Seifert fibre space \cite{BaiWoo1} and the factorisation in 
Proposition \ref{prop:fact} extends smoothly to $M^3$. 
\indent
When $n=3$ the $S^1$ action extends smoothly if $(M^4,g)$ is Einstein \cite{Pan-4to3} and, again, 
the factorisation extends smoothly to $M^4$. 
\end{rem}  

\indent
{}From Theorem \ref{thm:circle} and Proposition \ref{prop:fact} we obtain the following. 

\begin{thm} \label{thm:restr} 
Let  $\phi:(M^{n+1},g)\to(N^n,h)\:\,(n\geq3)$ be a non-constant harmonic morphism between 
compact oriented Riemannian manifolds.\\ 
\indent
Then all the Pontryagin numbers of $M^{n+1}$ are zero. In particular, the signature of $M^{n+1}$ 
is zero.\\
\indent 
Further, if $n\geq4$ then the Euler number of $M^{n+1}$ is zero. If $n=3$, then the Euler 
number of $M^4$ is even and is equal to the number of critical points of $\phi$. 
\end{thm} 
\begin{proof} 
If $n\geq4$ then, by Proposition \ref{prop:fact}, there exists a free $S^1$ action on $M^{n+1}$ 
whose orbits are connected components of the fibres of $\phi$. Hence, by the Hopf theorem the Euler 
number of $M^{n+1}$ is zero. Also, as is well-known (immediate consequence of \eqref{e:Bott}), 
all the Pontryagin numbers of $M^{n+1}$ are zero.\\
\indent
Suppose that $n=3$. If the set of critical points $\S$ is empty then the same argument as above 
implies that the Euler number and the Pontryagin number of $M^4$ are zero.\\
\indent
Suppose that $\S\neq\emptyset$ and let $x\in\S$ and $y=\phi(x)$. By Proposition \ref{prop:fact} 
we can assume that $\phi|_{M\setminus\phi^{-1}(\phi(\S))}$ has connected fibres. 
Let $B^3\subseteq N^3$ be a neighbourhood of $y\in N^3$ such that 
$B^3\cap\phi(\S)=\{y\}$ 
and which is diffeomorphic to the closed ball of radius two centred at zero in $\Re^3$. Then 
$\phi^{-1}(B^3)$ is a four-dimensional submanifold-with-boundary of $M^4$ such that 
$\phi^{-1}(B^3)\cap \S=\{x\}$. Furthermore, by Proposition \ref{prop:fact}, 
$\phi|_{\phi^{-1}(B^3)\setminus\{x\}}$ is the projection of an $S^1$-bundle over 
$B^3\setminus\{y\}$. Because $B^3$ is the cone over $S^2$, $\phi^{-1}(B^3)$ is the cone over its 
boundary. It easily follows (cf.\ \cite{Bai}) that the boundary of $\phi^{-1}(B^3)$ must be simply-connected.\\ 
\indent
Consider the embedding $S^2=\partial B^3\subseteq N^3\setminus\S$. Let $k\in\mathbb{Z}$ be the Chern number 
of the $S^1$ bundle $(\phi^{-1}(S^2),\,S^2,\,S^1)$. Then, if $k\neq0$, $\phi^{-1}(S^2)\cong S^3/\mathbb{Z}_{|k|}$ and, 
in particular, the fundamental group of $\phi^{-1}(S^2)$ is $\mathbb{Z}_{|k|}$. (If $k=0$ the bundle 
$(\phi^{-1}(S^2),\,S^2,\,S^1)$ is trivial.) But, $\phi^{-1}(S^2)$ 
is diffeomorphic to the boundary of $\phi^{-1}(B^3)$ which we have seen is simply-connected. 
It follows that $k=\pm1$ and, thus, we can suppose that $\phi|_{\phi^{-1}(B^3)\setminus\{x\}}$ is 
smoothly equivalent to the projection of the cylinder of the Hopf bundle $(S^3,S^2,S^1)$.\\ 
\indent
Thus, by taking, if necessary, the equivariant connected sum of $M$ and $-k\Co^2$, about $x\in M$ 
and $0\in\Co^2$, where $\Co^2$ is considered with its canonical circle action, we can suppose 
that on $\phi^{-1}(B^3)$ we have a smooth circle action having $x$ as a fixed 
point outside of which the action is free. By repeating this procedure about each point of $\S$ we 
obtain on $M^4$ a smooth circle action whose fixed point set is $\S$ outside which the action is free.\\ 
\indent 
By Theorem \ref{thm:circle}, the Pontryagin number of $M^4$ is zero and its Euler number is even and 
equal to the cardinal of $\S$.    
\end{proof} 

\begin{rem} 
With the same notations as in the proofs of Theorem \ref{thm:circle} and Theorem \ref{thm:restr} 
we have that $\epsilon(x)=k$ for each $x\in\S$.\\ 
\indent 
In fact, \eqref{e:sumzero} applied in the context of the proof 
of Theorem \ref{thm:restr} can be proved by applying Stokes' Theorem and the Chern-Weil Theorem, 
i.e. by assuming, if necessary, that $\phi$ has connected fibres and $V|_{M\setminus\S}$ has periodicity $2\pi$ then 
\begin{equation*}
0=\frac{1}{2\pi}\int_{N\setminus\,\cup_{x\in\S}\overset{\circ}{B}_x}\dif\!F
=\sum_{x\in\S} \Bigl( -\frac{1}{2\pi}\int_{\partial{}B_x}F \Bigr)=
\sum_{x\in\S}\epsilon(x) 
\end{equation*}
where $F\in\G(\Lambda^2(T^*(N\setminus\phi(\S)))$ is the curvature form of any principal 
connection on $\phi|_{M\setminus\S}$ and $B_x\subseteq N$ is a closed ball about each $\phi(x)$\,, $x\in\S$, 
such that $B_{x_1}\cap B_{x_2}=\emptyset$ for $x_1\neq x_2$. 
\end{rem} 

\indent
Let $P_g$ be the closed oriented surface of genus $g\geq0$ and let $n\geq2$. 
Since the Euler number of $\Co\!P^n$ is $n+1$, of 
$S^{2n}\times P_g$ is $4(1-g)$  and 
the signature of a $K3$ surface is $-16$ we have the following immediate consequence of 
Theorem \ref{thm:restr}.

\begin{cor} \label{cor:whatever} 
$K3$ surfaces, $\Co\!P^n$, $S^{2n}\times P_g$ ($n\geq2\,, g\neq1$) can never be the domain of a 
non-constant harmonic morphism with one-dimensional fibres whatever metrics we put on them. 
\end{cor} 

\begin{rem} 
1) Note that the projections $S^{2n+1}\times P_g\to\Co\!P^n\times P_g$ induced by the Hopf 
fibrations $S^{2n+1}\to\Co\!P^n$ are Riemannian submersions with 
totally-geodesic fibres, with respect to suitable multiples of their standard metrics, and are thus 
harmonic morphisms.\\
\indent
2) Other constructions of even-dimensional compact manifolds which cannot be the total space 
of a harmonic morphism with one-dimensional fibres whatever metrics we put on them can 
be easily obtained by using the product and/or the connected sum of manifolds. 
\end{rem}

\section{Harmonic morphisms with one-dimensional fibres\\ 
on compact Riemannian four-manifolds} 

\indent
The main result of this section is the following. 

\begin{thm} \label{thm:EinsteinASDscalar-flat} 
Let $(M^4,g)$ be a compact orientable Riemannian four-manifold and let 
$\phi:(M^4,g)\to(N^3,h)$ be a non-constant harmonic morphism.\\ 
\indent
Then the following assertions are equivalent\\
\indent
{\rm (i)} $(M^4,g)$ is half-conformally flat and its scalar curvature is zero,\\ 
\indent
{\rm (ii)} $(M^4,g)$ is Einstein and half-conformally flat,\\ 
\indent
{\rm (iii)} $(M^4,g)$ is Einstein and $\phi$ is submersive,\\ 
\indent
{\rm (iv)} $(M^4,g)$ is Ricci-flat.\\
\indent
Furthermore, if one of the assertions {\rm (i)}, {\rm (ii)}, {\rm (iii)} or {\rm (iv)} holds then, 
up to homotheties and Riemannian coverings, $\phi$ is the canonical projection $T^4\to T^3$ 
between flat tori. 
\end{thm} 
\begin{proof} 
Choose one of the orientations of $M^4$ and let $\omega_g$ be the corresponding volume-form 
with respect to $g$.\\ 
\indent
Let $p_1[M]$ be the Pontryagin number of $M^4$. By the Chern-Weil Theorem we have 
\begin{equation} \label{e:p1} 
p_1[M]=\frac{1}{4\pi^2}\int_M\big(\,|W^+|^2-|W^-|^2\,\big)\,\omega_g
\end{equation} 
where $W$ is the Weyl tensor of $(M^4,g)$ and $W^+$\,, $W^-$ are its self-dual and anti-self-dual 
components, respectively (see \cite[13.8]{Bes}).\\ 
\indent
By Theorem \ref{thm:restr}, $p_1[M]=0$ and hence $W^{\pm}=0\iff W^{\mp}=0$. Thus, if $(M^4,g)$ is 
half-conformally flat then it is conformally flat.\\ 
\indent
Now, recall that, by the Gauss-Bonnet Theorem, the Euler number of $M^4$ is given by (see \cite{Bes}, \cite{LeBr1}): 
\begin{equation} \label{e:Euler} 
\chi[M]=\frac{1}{8\pi^2}\int_M\Big(\,\frac{s^2}{24}-\frac{|\RicM_{\,0}|^2}{2}+|W^+|^2
+|W^-|^2\,\Big)\,\omega_g
\end{equation} 
where $s$ is the scalar curvature of $(M^4,g)$ and $\RicM_{\,0}$ is the trace-free part of the Ricci tensor 
of $(M^4,g)$.\\
\indent
If $(M^4,g)$ is half-conformal flat and its scalar curvature is zero then \eqref{e:Euler} becomes 
\begin{equation*} 
\chi[M]=-\frac{1}{16\pi^2}\int_M|\RicM_{\,0}|^2\,\omega_g\;.
\end{equation*} 
But, by Theorem \ref{thm:restr}, $\chi[M]\geq0$ and hence $(M^4,g)$ must be Einstein.\\   
\indent
If $(M^4,g)$ is 
Einstein and half-conformally flat then $(M^4,g)$ has constant sectional curvature $k^M$ (see \cite{Bes}). 
By \cite[Proposition 3.3(ii)]{Pan} and \cite[Proposition 3.6]{Pan-4to3} we cannot have $k^M<0$. If $k^M>0$ 
then, up to homotheties, the universal cover of $(M^4,g)$ is $S^4$, a situation which cannot occur 
(see \cite[Section 3]{Bry}). Hence $(M^4,g)$ must be flat.\\ 
\indent
If $(M^4,g)$ is Einstein and $\phi$ is submersive then, by Theorem \ref{thm:restr}, the Euler number of 
$M^4$ is zero. Thus, \eqref{e:Euler} 
implies that $(M^4,g)$ is flat.\\  
\indent
If $(M^4,g)$ is Ricci-flat then, as a consequence of \cite[Proposition 3.6]{Pan-4to3}, $\phi$ must be 
submersive, since on a compact Ricci-flat  manifold any Killing vector field is parallel (see \cite{Kob}). 
Then, by Theorem \ref{thm:restr}, $\chi[M]=0$. Now, \eqref{e:Euler} implies that $(M^4,g)$ is flat. 
The last assertion follows from \cite[Theorem 3.4]{Pan} 
and an argument as in the proof of \cite[Theorem 3.8]{Pan-4to3}. 
\end{proof} 

\begin{rem} 
The fact that assertion (iii) implies that, up to homotheties and Riemannian coverings, $\phi$ is 
the canonical projection $T^4\to T^3$ between flat tori is the result of \cite[Theorem 3.8]{Pan-4to3}. 
\end{rem} 

\indent
{}From the proof of Theorem \ref{thm:EinsteinASDscalar-flat} we obtain the following. 

\begin{prop}
Let $\phi:(M^4,g)\to(N^3,h)$ be a non-constant harmonic morphism defined on an orientable compact 
Riemannian four-manifold.\\ 
\indent
Then $(M^4,g)$ is half-conformally flat if and only if it is conformally flat. 
\end{prop}

\indent
Next we give a sufficient condition for the total space of a harmonic morphism $\phi:(M^4,g)\to(N^3,h)$ 
to be half-conformally flat. Note that this does not require any compactness or completeness 
assumption. 

\begin{prop} \label{prop:ASD} 
Let $(M^4,g)$ be an orientable Einstein four-manifold.\\ 
\indent
Suppose that there exists  a non-constant 
harmonic morphism $\phi:(M^4,g)\to(N^3,h)$ to an orientable Einstein three-manifold $(N^3,h)$.\\ 
\indent
Then $(M^4,g)$ is half-conformally flat. 
\end{prop} 
\begin{proof} 
Let $x\in M$ be a regular point of $\phi$ and $y=\phi(x)$. Let $Y_0\in T_yN$ be a unit vector. 
Because $(N^3,h)$ is of constant curvature, by \cite{BaiWoo-spfo},  
there exists an open neighbourhood $U$ of $y$ and a submersive 
harmonic morphism $\psi_{Y_0}:(U,h|_U)\to P^2$ with values in some Riemann surface $P^2$ such that its 
fibre through $y$ is tangent to $Y_0$. Then for any other unit vector $Y\in T_yN$   
we can compose $\psi_{Y_0}$ with an isometry to obtain a submersive harmonic morphism 
$\psi_Y$ whose fibre through $y$ is tangent to $Y$.\\ 
\indent
Then, $\psi_Y\circ\phi:(\phi^{-1}(U),g|_{\phi^{-1}(U)})\to P^2$ is a submersive harmonic morphism from an 
orientable Einstein four-manifold to a Riemann surface. By \cite[Theorem 1.1]{Woo-4d}, there exists 
an (integrable) Hermitian structure $J_Y$ on $(\phi^{-1}(U),g|_{\phi^{-1}(U)})$ with respect to which  
$\psi_Y\circ\phi:(\phi^{-1}(U),J_Y)\to P^2$ is holomorphic.  By restricting, if necessary, the family of 
$Y$ to an open subset of the unit 
sphere in $(T_yN,h_y)$ we can suppose that all the $J_Y$ induce the same orientation $\sigma$ on 
$\phi^{-1}(U)$. 
By the Riemannian Goldberg-Sachs Theorem (see \cite{ApoGau}) $W^+$ is degenerate.  
Now, by a result of A.~Derdzi\'nski \cite{Der} either $W^+=0$ or there is just exactly one pair $\pm J$ 
of (oriented) complex structures compatible with $g|_{\phi^{-1}(U)}$. But the latter cannot occur because,  
 obviously, if $Y_1\neq\pm Y_2$ then $J_{Y_1}\neq\pm J_{Y_2}$. Hence $W^+=0$ and the proof follows. 
\end{proof} 

\indent
By combining Theorem \ref{thm:EinsteinASDscalar-flat} and Proposition \ref{prop:ASD} we immediately obtain a new proof 
for the following result from \cite[Theorem 4.11]{Pan-4to3}. 

\begin{thm} \label{thm:compactbothEinstein}
Let $(M^4,g)$ be a compact Einstein four-manifold and let $(N^3,h)$ be a Riemannian three-manifold
with constant curvature.\\ 
\indent
Let $\phi:(M^4,g)\to(N^3,h)$ be a non-constant harmonic morphism.\\ 
\indent
Then, up to homotheties and Riemannian 
coverings, $\phi$ is the canonical projection $T^4\to T^3$ between flat tori.\\
\indent
In particular there exists no harmonic morphism with one-dimensional 
fibres from a compact Einstein manifold of dimension four to $(S^3,{\rm can})$\,.
\end{thm} 

\begin{rem} 
There exists a harmonic morphism from $S^4$ considered with a metric which is conformally equivalent 
to the canonical one to $(S^3,{\rm can})$ (see \cite{BaiRat}).  Note that, in \cite{BaiRat}, 
there is also constructed a harmonic morphism from $(S^4,{\rm can})$ to $(S^3,h)$ where $h$ is continuous 
and smooth away from the poles. 
\end{rem}

\section{Non-existence of harmonic morphisms on\\*
compact Hermitian-Einstein four-manifolds} 

In this section we prove the following. 

\begin{thm} \label{thm:Hermitian-Einstein}
Let $(M^4,g,J)$ be a compact Hermitian-Einstein four-manifold 
and let $\phi:(M^4,g)\to(N^3,h)$ be a non-constant harmonic
morphism to an arbitrary Riemannian manifold.\\ 
\indent
Then, up to homotheties and Riemannian 
coverings, $\phi$ is the canonical projection $T^4\to T^3$ between flat tori. 
\end{thm}
\begin{proof} 
If $\phi$ is submersive then the proof follows from \cite[Theorem 3.8]{Pan-4to3}.\\ 
\indent
Suppose that $\phi$ has critical points. Then, by \cite[Proposition 3.6]{Pan-4to3}, 
$(M^4,g)$ has positive scalar curvature and there exists a \emph{smooth} Killing vector field $V$ 
tangent to the fibres of $\phi$ whose zero set is equal to the set of critical 
points of $\phi$.\\ 
\indent
Suppose that $(M^4,g,J)$ is K\"ahler-Einstein. Its first Chern class is positive. 
Then, by \cite[Proposition 3.2]{Hit} or \cite[Theorem I.1.3.4]{Yau}, $(M^4,J)$ is either 
$\Co\!P^1\times\Co\!P^1$ or is obtained from $\Co\!P^2$ by blowing-up $r$ distinct points, 
$0\leq r\leq8$; such a blow-up has signature $1-r$.  
 But, by Theorem \ref{thm:restr}, $M^4$ has signature zero and hence either 
$M^4$ is $\Co\!P^1\times\Co\!P^1$ or is obtained from $\Co\!P^2$ by blowing-up 
one point. In the latter case, $(M^4,J)$ is biholomorphic to 
$\Co\!P^2$\!\#$\overline{\Co\!P^2}$; but this admits \emph{no} K\"ahler-Einstein metric 
(see \cite[11.54, 11.56]{Bes}). If $(M^4,J)=\Co\!P^1\times\Co\!P^1$ then, 
by \cite[Proposition I.1.4.5]{Yau}, $(M^4,g,J)$ is homothetic to $\Co\!P^1\times\Co\!P^1$.\\ 
\indent
Hence we may suppose that $(M^4,g,J)$ is isometric to $\Co\!P^1\times\Co\!P^1$.  
Let $(x,y)\in\Co\!P^1\times\Co\!P^1$ be a critical point of $\phi$ (note that by 
Theorem \ref{thm:restr}, in this case, $\phi$ must have exactly four critical points). 
But $(x,y)$ is also an isolated zero of  $V$. Because $V$ is a Killing vector field,  
$(\nabla V)_{(x,y)}$ preserves each of the summands in the orthogonal decomposition 
$T_{(x,y)}(\Co\!P^1\times\Co\!P^1)=T_x\Co\!P^1\oplus T_y\Co\!P^1$. Now, recall that 
$(\nabla V)_{(x,y)}$ is an orthogonal complex structure on $T_{(x,y)}(\Co\!P^1\times\Co\!P^1)$.   
Hence we can suppose that $(\nabla V)_{(x,y)}=J_{(x,y)}$. It follows that, by composing 
$\phi$ with the inverse of the map $\Co\!P^1\times\Co\!P^1\setminus\{(-x,-y)\}\to\Co\times\Co$,  
given by stereographic projection on each factor, we get a harmonic morphism 
$\psi:(\Co^2,\bar{g})\to(N^3,h)$ 
where $$\bar{g}=\frac{4}{(1+|z|^2)^2}|\dif\!z|^2+\frac{4}{(1+|w|^2)^2}|\dif\!w|^2\;.$$ 
Since the stereographic projection is conformal $\psi$ is induced by the 
canonical circle action on $\Co^2$; note that this is an isometric action with respect to $\bar{g}$. 
Let $\bar{V}$ be its infinitesimal generator and let $\l$ be the 
dilation of $\psi$. Then, up to a multiplicative constant we have $\bar{g}(\bar{V},\bar{V})=\l^2$ 
(see \cite{Bry}, \cite[Section 2]{Pan}) and $\psi^*(h)$ must be equal to the horizontal component of 
$\l^2\bar{g}$. However, it can be checked without difficulty that then $h$ \emph{cannot} be extended 
over $\psi(0)$.\\ 
\indent
We have thus proved that $(M^4,g,J)$ cannot be K\"ahler-Einstein. Because $M^4$ has signature zero, 
from the main result of \cite{LeBr} it follows that $(M^4,g)$ is $\Co\!P^2$ with one point 
blown-up endowed with the Page metric. However, from the discusion in \cite{GibHaw} it follows 
that none of the Killing vector fields of the Page metric has isolated fixed points which 
are isolated as singularities as well so that this case is not possible either. 
\end{proof} 

\begin{rem} Thus there is no harmonic morphism from $\Co\!P^2\#-\Co\!P^2$ endowed with the 
Page metric to any $3$-manifold; there is, however a harmonic morphism from $\Co\!P^2\#-\Co\!P^2$ endowed with the Page metric to $\Co\!P^1$ (see \cite{BaiWoo2}). 
\end{rem}

\end{document}